\documentclass[12pt]{article}
\usepackage{amsmath,amsthm}
\usepackage{amssymb,latexsym}
\usepackage{enumerate}

\theoremstyle{definition}

\numberwithin{equation}{section}

\newcounter{theoremcounter}
\newenvironment{theorem}{\refstepcounter{theoremcounter}{\vspace*{3mm}
   \par\noindent \bf Theorem \thetheoremcounter.}\it}{}

\newenvironment{corollary}{\refstepcounter{theoremcounter}{\vspace*{3mm}
   \par\noindent \bf Corollary \thetheoremcounter.}\it}{}

\renewcommand{\=}{\overset{ \text{def} }=}
\DeclareMathOperator{\ddh}{dh}
\newcommand{\thmref}[1]{Theorem~\ref{#1}}

\renewenvironment{proof}[1][Proof]{\begin{trivlist}
\item[\hskip \labelsep {\bfseries #1}]}{\end{trivlist}}
\renewcommand{\qed}{\nobreak \ifvmode \relax \else
      \ifdim\lastskip<1.5em \hskip-\lastskip
      \hskip1.5em plus0em minus0.5em \fi \nobreak
      \vrule height0.75em width0.5em depth0.25em\fi}

\newcommand{\R}{{\mathbf{R}}}

\newcommand{\E}{{\mathbf{E}}}
\newcommand{\F}{{\mathcal{F}}}

\begin{document}


\baselineskip=17pt


\title{Extremal Lipschitz functions in the deviation inequalities from the mean}

\author{Dainius Dzindzalieta$^1$}

\maketitle


\renewcommand{\thefootnote}{}

\footnote{$^1$ Institute of Mathematics and Informatics, Vilnius University,
Vilnius, Lithuania.
E-mail: dainiusda@gmail.com}
\footnote{$^1$ This research was funded by a grant (No. MIP-12090) from the Research Council of Lithuania}


\renewcommand{\thefootnote}{\arabic{footnote}}
\setcounter{footnote}{0}


\begin{abstract}
We obtain an optimal deviation from the mean upper bound
\begin{equation}
D(x)\=\sup_{f\in \F}\mu\{f-\E_{\mu} f\geq x\},\qquad\ \text{for}\ x\in\R\label{abstr}
\end{equation} where $\F$ is the class of the integrable, Lipschitz functions on probability metric (product) spaces.
As corollaries we get exact solutions of $\eqref{abstr}$ for Euclidean unit sphere $S^{n-1}$ with a geodesic distance and a normalized Haar measure, for $\R^n$ equipped with a Gaussian measure and for the multidimensional cube, rectangle, torus or Diamond graph equipped with uniform measure and Hamming distance. We also prove that in general probability metric spaces the $\sup$ in $\eqref{abstr}$ is achieved on a family of distance functions.
\end{abstract}

\section{Introduction}

\label{pirmassk}

Let us recall a well known result for Lipschitz functions on probability metric spaces.
We denote by $\F=\F(V)$ the class of integrable, i.e. $f\in L_1(V,d,\mu)$, $1$-Lipschitz functions $f:V\rightarrow\R$, i.e., such that $~{|f(u)-f(v)|\leq d(u,v)}$ for all $u,v\in V$. Let $M_f$ be a median of the function $f$, i.e., a number such that $\mu\{f(x)\leq M_f\}\geq \frac12$ and $~{\mu\{x:f(x) \geq M_f\}\geq \frac12}$. Given probability metric space $(V,d,\mu)$
$$
\sup_{f\in \F}\mu\{f-M_f\geq x\}\qquad \text{for}\ x\in\R
$$ is achieved on a family, say $\F^{\ast}$, of distance functions $f(u)=-d(A,u)$ with measurable $A\subset V$ (for a book type exposition of results we refer reader to \cite{ledoux2011probability, milman2002asymptotic}). From this it is easy to deduce that this problem is equivalent to the following vertex isoperimetric problem. Given $t\geq0$ and $h\geq 0$,
\begin{equation}
\text{minimize}\ \mu(A^h)\ \text{over all}\ A\subset V\ \text{with}\ \mu(A)\geq t, \label{ext}
\end{equation} where $A^h=\{u\in V: d(u,A)\leq h\}$ is an $h$-enlargement of $A$.

Following  $\cite{bezrukov2002local}$ we say that a space $(V,d,\mu)$ is \textit{isoperimetric} if for every $t\geq0$ there exists a solution, say $A_{\text{opt}}$, of $\eqref{ext}$ which does not depend on $h$.

\bigskip

However, as was pointed out by Talagrand \cite{talagrand1996new} in practise it is easier to deal with expectation $\E f$ rather than $M_f$.
In order to get results for the mean instead of median two different techniques are usually used. One way is to evaluate the distance between median and mean, another is to use a martingale technique (see \cite{mcdiarmid1989method, talagrand1996new, ledoux1996talagrand,bentkus2007measure} for more detailed exposition of the results). Unfortunately, none of them lead to tight bounds for the mean.

\bigskip

In this paper we find tight deviation from the mean bounds
\begin{equation}
D(x)\=\sup_{f\in \F}\mu\{f-\E_{\mu} f\geq x\},\qquad\ \text{for}\ x\in\R\label{bendra}
\end{equation} for the class $\F$ of all Lipschitz functions on a probability metric space $(V,d,\mu)$. Here a probability metric space means that a measure $\mu$ is Borel and normalized, $\mu(V)=1$. If we change $f$ to $-f$ we get that
$$
D(x)=\sup_{f\in \F}\mu\{f-\E_{\mu} f\leq -x\}\qquad \text{for}\ x\in\R.
$$
 Note that the function $D(x)$ depends also on $(V,d,\mu)$.

\bigskip

We first state a general result for all probability metric spaces. The $\sup$ in $\eqref{bendra}$ is achieved on a family, say $\F^{\ast}$, of distance functions $f(u)=-d(A,u)$ with measurable $A\subset V$. Note that $\F^{\ast}\subset\F$.
\begin{theorem}\label{general}
Given a probability metric space $(V,d,\mu)$, we have
$$ \sup_{f\in\F}\mu\{f-\E_{\mu} f\geq x\}=\sup_{f\in\F^{\ast}}\mu\{f-\E_{\mu} f\geq x\}\qquad\ x\in\R.
$$
\end{theorem}
\begin{proof}
Fix $x\in\R$.  Let $f\in\F$ and $B=\{u:f(u)-\E_{\mu} f\geq x\}$. Without loss of generality we suppose that $f\geq 0$  for $u\in B$ and $f(u)=0$ for at least one $u\in B$. Let $g(u)=-d(B,u)$. Since $f$ is Lipschitz function,
we get that $f(u)\geq g(u)$ for all $u\in V$, so $\E_{\mu} f\geq \E_{\mu} g$. From this we get that $B\subset\{u: g(u)-\E_{\mu} g\geq x\}$, thus
 $
 \mu\{f-\E_{\mu} f\geq x\} \leq \mu\{g-\E_{\mu} g\geq x\}.
 $
Since $f$ is arbitrary the statement of $\thmref{general}$ follows.\qed
\end{proof}

In the special case when  $V=\R^n$ and $\mu=\gamma_n$ is a standard Gaussian measure, $\thmref{general}$ was proved by Bobkov $\cite{bobkov2003localization}$.

\bigskip

Our main result is the following theorem.
\begin{theorem}\label{pagr}
 For all isoperimetric probability metric spaces $(V,d,\mu)$ we have
$$
D(x)=\mu\{f^{\ast}-\E_{\mu} f^{\ast}\geq x\}\qquad \text{for}\ x\in\R,\notag
$$
 where $f^{\ast}(u)=-d(A_{\text{opt}},u)$ is a distance function from some extremal set $A_{\text{opt}}$. It turns out that $\mu(A_{\text{opt}})=D(x)$.
\end{theorem}

\bigskip

\begin{proof}
We have
\begin{eqnarray*}
\E_{\mu} d(A,\cdot)&=&\int_0^\infty \left(1 - \mu\left\{A^h\right\}\right) \ddh\\
&\leq& \int_0^\infty \left(1-\mu\left\{A_{\text{opt}}^h\right\}\right) \ddh=\E_{\mu} d(A_{\text{opt}},\cdot).
\end{eqnarray*}
Thus by $\thmref{general}$
$$
\sup_{f\in \F}\mu\{f-\E_{\mu} f\geq x\}=\sup_{f\in \F^{\ast}}\mu\{f-\E_{\mu} f\geq x\}= \mu\{f^{\ast}-\E_{\mu} f^{\ast}\geq x\}
$$ with $f^{\ast}(u)=-d(A_{\text{opt}},u)$. Note that the second equality holds since looking for the $\sup$ we can consider only the functions $f^{\prime}$ from the class $\F^{\ast}$ such that $\mu\{x:f^{\prime}(x)=0\}=D_n(x)=\mu\{x:f^{\ast}(x)=0\}$.
Which completes the proof of $\thmref{pagr}$.\qed
\end{proof}

\section{Isoperimetric spaces and corollaries}
In this section a short overview of results on the vertex isoperimetric problem described by \eqref{ext}.
We also state a number of corollaries following from $\thmref{general}$ and $\thmref{pagr}$.

A typical and basic example of isoperimetric spaces is the Euclidean unit sphere $S^{n-1}=\{x\in\R^n:\ \sum_{i=1}^n |x_i|^2=1\}$ with a geodesic distance $\rho$ and a normalized Haar measure $\sigma_{n-1}$. P. L{\'e}vy $\cite{levy1951}$ and E. Schmidt $\cite{schmidt48}$ showed that if $A$ is a Borel set in $S^{n-1}$ and $H$ is a cap (ball for geodesic distance $\rho$) with the same Haar measure $~{\sigma_{n-1}(H)=\sigma_{n-1}(A)}$, then
\begin{equation}\label{levy}
\sigma_{n-1}(A^h)\geq \sigma_{n-1}(H^h)\qquad \text{ for\ all }\ h>0.
\end{equation}
Thus $A_{\text{opt}}$ for the space $(S^{n-1},\rho,\sigma_{n-1})$ is a cap. We refer readers for a short proof  of $\eqref{levy}$ to $\cite{beny, flm1977}$. The extension to Riemannian manifolds  with strictly positive curvature can be found in $\cite{grom1}$. Note that if $H$ is a cap, then $H^h$ is also a cap, so we have an immediate corollary.
\begin{corollary}
 For a unit sphere $S^{n-1}$ equipped with normalized Haar measure $\sigma_{n-1}$ and geodesic distance we have
$$
  D(x)=\sigma_{n-1}\{f^{\ast}-\E_{\mu} f^{\ast}\geq x\}\qquad \text{for}\ x\in\R,
$$ where  $f^{\ast}(u)=-d(A_{\text{opt}},u)$ and $A_{\text{opt}}$ is a cap.
\end{corollary}

Probably the most simple non-trivial isoperimetric space is $n$-dimensional discrete cube $~{C_n=\left\{0,1\right\}^n}$ equipped with uniform measure, say $\mu$, and Hamming distance. Harper $\cite{harper1964optimal}$ proved that some number of the first elements of $C_n$ in the simplical order is a solution of $\eqref{ext}$. Bollobas and Leader $\cite{bollobas1991isoperimetric}$ extended this result to multidimensional rectangle. Karachanjan and Riordan \cite{karachanjan1982discrete, riordan1998ordering} solved the problem \eqref{ext} for multidimensional torus. Bezrukov considered
 powers of the Diamond graph $\cite{bezrukov2008vertex}$ and powers of cross-sections $\cite{bezrukov2002local}$.  We state the results for discrete spaces as corollary.
\begin{corollary}
 For  discrete multidimensional cube, rectangle, torus and Diamond graph equipped with uniform measure and Hamming distance we have
 $$
D(x)=\mu\{f^{\ast}-\E_{\mu} f^{\ast}\geq x\}\qquad \text{for}\ x\in \R,
$$
 where $f^{\ast}(u)=-d(A_{\text{opt}},u)$ and $A_{\text{opt}}$ are the sets of some first elements in corresponding orders. In particular, for $n$-dimensional discrete cube with Hamming distance $A_{\text{opt}}$ is a set of some first elements of $C_n$ in simplical order.
\end{corollary}

There is a vast of papers dedicated to bound $D(x)$ for various discrete spaces. We mention only $\cite{bobkov1997isoperimetric,ledoux2001concentration,talagrand1995concentration}$ among others. In $\cite{bezrukov2008vertex, ledoux2011probability}$ a nice overview of isoperimetric spaces and bounds for $D(x)$ are provided.

\bigskip
Another important example of isoperimetric spaces comes from Gaussian isoperimetric problem. Sudakov and Tsirel'son $\cite{sudakov1978extremal}$ and Borell $\cite{borell1975brunn}$ discovered that if $\gamma_n$ is a standard Gaussian measure on $\R^n$ with a usual Euclidean distance function $d$, then $(\R^n,d,\gamma_n)$ is isoperimetric. In $\cite{sudakov1978extremal}$ and $\cite{borell1975brunn}$ it was shown that among all subsets $A$ of $\R^n$ with  $t\geq\gamma_n(A)$, the minimal value of $\gamma_n(A^h)$ is attained for half-spaces of measure $t$. Thus we have the following corollary of $\thmref{general}$ and $\thmref{pagr}$.

\begin{corollary} For a Gaussian space  $(\R^n,d,\gamma_n)$ we have
$$
D(x)= \gamma_n\{f^{\ast}-\E_{\gamma_n} f^{\ast} \geq x\}\qquad\text{for}\ x\in \R, \label{bobkovo}
$$ where $~{f^{\ast}(u)=-d(A_{\text{opt}},u)}$ is a distance function from a half-space of space $\R^n$.
\end{corollary}

The latter result was firstly proved by Bobkov $\cite{bobkov2003localization}$. We also refer for further investigations of extremal sets on $\R^n$ for some classes of measures to $\cite{barthe2007isoperimetry,bobkov1997b}$ among others.

\bibliographystyle{alpha}

\end{document}